\documentclass[review]{elsarticle}%
\usepackage{graphicx}
\usepackage{amsmath}
\usepackage{amssymb}
\usepackage{colortbl}
\usepackage{amsthm}
\usepackage{arydshln}
\usepackage{booktabs}
\usepackage{lineno,hyperref}
\modulolinenumbers[5]

\usepackage{setspace}

\newcommand{\diag}{\mathop{\textrm{diag}}}

\newtheorem{thm}{Theorem}

\newtheorem{rem}{Remark}

\newcommand{\ot}{\overline{\tau}}

\newcommand{\clr}{\color{black} }

\journal{}

\begin{document}

\begin{frontmatter}

\title{{\clr{LMI approach to global stability analysis of stochastic delayed Lotka-Volterra models}}}
\author[kk]{Krisztina Kiss\corref{mycorrespondingauthor}}
\cortext[mycorrespondingauthor]{Corresponding author}
\ead{kk@math.bme.hu}
\author[kk]{\'Eva Gyurkovics}
\ead{gye@math.bme.hu}

\address[kk]{Mathematical Institute, Budapest University of Technology and Economics, Budapest, Pf. 91,
1521, Hungary}
\begin{abstract}This paper is devoted to the stability analysis of {\clr{an} $n$} species Lotka-Volterra {\clr{system}}
with discrete and distributed delays. Stochastic perturbations  to the parameters of the model are allowed. Sufficient conditions for the almost sure global asymptotic stability of the 
{\clr {positive}} equilibrium are derived in terms of LMIs. The efficiency of the proposed method is illustrated by numerical examples.
\end{abstract}

\end{frontmatter}

\section{Introduction}
In the past decades, one of the
most popular models in mathematical biology has been the Lotka-Volterra model
that have been studied in a huge number of works. In particular, the books \cite{maokonyv}-
\cite{smithkonyv} are good references in this area.
A large class of models is given by
ordinary differential equations,
  but it is often more realistic to use delayed functional
differential equations (FDEs) to describe such models (e.g. \cite{Chen2016}-
-\cite{VargasAbs2015}).
Consider the n-species Lotka-Volterra model of the form
\begin{align}\label{VL1}
 \dot{u}_i (t)  = u_i(t)\left[\rho_i - \sum_{j=1}^n a_{ij}u_j(t)-\sum_{j=1}^{n} a_{ij}^d u_j(t-\tau_{ij}(t)) -\sum_{j=1}^{n} a_{ij}^D \int_{t-\overline{\tau}_{ij}}^{t} e^{\alpha_{ij}(\eta-t)}u_j(\eta)d\eta\right],
  \hspace{0.3cm} i=1,\ldots, n
\end{align}
consisting both discrete 
and distributed delays. Here $u_i(t), \; (i=1,\ldots,n)$ represent the population sizes,
the parameters $a_{ij}, a_{ij}^d, a_{ij}^D$ are the so-called interaction coefficients, {\clr{
 $\rho_i>0$ is the $i$th intrinsic growth rate,
 while the values $\alpha_{ij} \geq 0$ play the role of the weighting parameters of the distributed delays.}}
The discrete delays $\tau_{ij}$ are supposed to be differentiable functions of the time satisfying conditions
\begin{equation}\label{tau1}
  0< \tau_{ij}(t) \leq \overline{\tau}_{ij}, \hspace{0.3cm} \dot{\tau}_{ij}(t)\leq \overline{\tau}_{ij}^d \hspace{0.3cm}
\end{equation}
with known constants $\overline{\tau}_{ij}, \, \overline{\tau}_{ij}^d${\clr{.}}
 {\clr{Define $\overline{\tau}= \max_{1\leq i,j \leq n} \overline{\tau}_{ij}.$}}

Let us suppose that \eqref{VL1} has a positive equilibrium state $u^*=(u_1^*, u_2^*,\dots, u_n^*)^T \in \mathbf{R}_{+}^{n},$ i.e.
   $\tilde{A} u^{\ast}  = \rho,$
where the notations
 $\tilde{A}=A+A^d+A_{\beta}^D,$ 
 $A =(a_{ij})_{n\times n},$ 
 $A^d=(a_{ij}^d)_{n\times n},$ 
 $A_{\beta}^D=(\beta_{ij}a_{ij}^D)_{n\times n},$ 
 $ \rho=(\rho _1,\ldots,\rho _n)^T, $ 
with
\begin{equation}
  \beta_{ij}= \int_{t-\overline{\tau}_{ij}}^{t} e^{\alpha_{ij}(\eta-t)}d\eta=
  \left\{ \begin{array}{ll}
   \frac{1}{\alpha_{ij}}(1-e^{-\alpha_{ij}\overline{\tau}_{ij}}),
   & \mbox{if $\alpha_{ij} \neq 0$,} \\
   \overline{\tau}_{ij},
   & \mbox{if $\alpha_{ij}=0$}
                     \end{array}
          \right.
 \label{beta}
\end{equation}
have been used.

The population systems are almost always subjected to environmental noises (see e.g. \cite{maokonyv}, \cite{LiuPhys2015}-\cite{XiongAML2019}). Similarly to \cite{mao2005}, we will take into account random fluctuations, namely white noise, affecting on $\rho_i$, and depending on how much the current population sizes differ from the equilibrium state. Thus, we replace $\rho_i$ by $\rho_i + \sum_{j=1}^{n}\sigma_{ij}(u_j(t)-u_j^*)\dot{w}_i(t)$, where $\sigma^2_{ij}$ denotes the intensity of the noise and  $w(t)=(w_1(t), w_2(t), \dots, w_n(t))^T$ is a Brownian motion defined on a complete probability space $(\Omega, \mathcal{F}, \cal{P})$ with a filtration $\left\{ \mathcal{F}_t \right\}_{t\geq 0}$ satisfying the "usual conditions" (see e.g. in \cite{LiuPhys2015}). In this way, we obtain the stochastic system, which can be written as
\begin{align}\label{VL2}
 d{u}_i (t)
 =& u_i(t)\left\{ \left[\rho_i -  \sum_{j=1}^n a_{ij}u_j(t)-\sum_{j=1}^{n} a_{ij}^d u_j(t-\tau_{ij}(t))
 - \sum_{j=1}^{n} a_{ij}^D \int_{t-\overline{\tau}_{ij}}^{t} e^{\alpha_{ij}(\eta-t)}u_j(\eta)d\eta \right]dt \right. \nonumber \\
&\hspace{3.9cm} \left. +\sum_{j=1}^{n}\sigma_{ij}(u_j(t)-u_j^*)d{w}_i(t)
\right\},
 \hspace{2.5cm} i=1, \ldots, n.
\end{align}
In order to write model \eqref{VL2}
in a more compact form, we introduce some notations. Let $u^d(t), \; u^D(t)\in \mathbf{R}^{n^{2}}$ be vectors having elements
\begin{equation}\label{Gjel1}
  u_k^d(t)=u_j(t-\tau_{ij}(t)),\hspace{0.3cm}
  u_k^D(t)=\int_{t-\overline{\tau}_{ij}}^{t} e^{\alpha_{ij}(\eta-t)}u_j(\eta)d\eta, \hspace{0.5cm} \mbox{for}\hspace{0.2cm} k=(i-1)n+j, \; i,j=1,\ldots,n.
\end{equation}
Let the matrices $\mathcal{A}^d, \mathcal{A}^D\in \mathbf{R}^{n\times n^{2}}$ be defined  by
 $ \mathcal{A}^d=\diag\{{\underline{a}_1^d},\ldots,{\underline{a}_n^d}\}$ 
 and 
 $\mathcal{A}^D=\diag\{{\underline{a}_1^D},\ldots,{\underline{a}_n^D}\},$
where $\underline{a}_i^d$ and $\underline{a}_i^D$ are the $i$th row vectors of the matrices $A^d, \, A^D,$ respectively. Further, let $g : \mathbf{R}^n \rightarrow \mathbf{R}^{n \times n}$ be a diagonal matrix with $g_{ii}(u)=\sum_{j=1}^{n}\sigma_{ij}u_j$, and zero otherwise. Then system \eqref{VL2} can be written as  
\begin{align}\label{VLm}
 du(t)  = \diag\{u_1(t), \ldots,u_n(t)\}\left\{\left[ \rho -{A}u(t)-{\mathcal{A}}^d{u}^d(t)-{\mathcal{A}^D}{u}^D(t))\right]dt+g(u(t)-u^{\ast})dw(t)\right\}.
\end{align}
We shall consider system \eqref{VLm} with the initial condition $u(t)=\varphi_0(t),$ if $t\in [-\overline{\tau},0]$ and $\varphi_0 \in C([-\overline{\tau},0],\mathbf{R}_{+}^{n}).$
\emph{Our aim is to derive sufficient conditions for ensuring that
\begin{itemize}
  \item equation \eqref{VLm} has global positive solution almost surely for any initial function $\varphi_0 \in C([-\overline{\tau},0],\mathbf{R}_{+}^{n});$
  \item the equilibrium state of \eqref{VLm} is {\clr{almost surely}} globally asymptotically stable in  $\mathbf{R}_{+}^{n}$.
\end{itemize}
}
Based on some new development{\clr{s}} in the field, a new Lyapunov-Krasovskii functional is introduced, and the stability condition is given in terms of linear matrix inequalities (LMIs). To the best of our knowledge,  only variational system based local results have been given up to now by means of LMIs in the literature (see \cite{DongJIA2019}, and the references therein). 
The result obtained in this work demonstrates that LMI can be applied for investigating the stability behaviour of the \emph{nonlinear} Lotka-Volterra equation.
\section{Main results}
We shall first formulate a condition under which system \eqref{VLm} has a unique global positive solution a.s.
\begin{thm} \label{thm:1}
Assume that the discrete delay functions satisfy \eqref{tau1} and the condition
\begin{equation}\label{l-th1}
 \overline{\tau}^d := \max_{1\leq i,j \leq n} \overline{\tau}_{ij}^d <1,
\end{equation}
and the noise intensity matrix $\sigma=(\sigma_{ij})_{n\times n}$ has the property that $\sigma_{ii}>0$ for $i=1,\ldots,n,$ and $\sigma_{ij}\geq 0$ for $i,j=1,\ldots,n.$
Then, for any initial function $\varphi_0 \in C([-\overline{\tau},0],\mathbf{R}_{+}^{n})$ equation \eqref{VLm} has a unique positive solution $u(t) $ on $[-\overline{\tau},\infty), $ and the solution remains in $\mathbf{R}_{+}^{n}$ with probability 1.
\end{thm}
{\bf Proof.} The proof follows the same line as \cite{mao2005} and \cite{LiuPhys2015},
therefore the details are omitted to save space.  $\square $

Now we turn to the problem of stochastic asymptotic stability of the equilibrium state of \eqref{VLm}.
To this end,
shifting the origin to the equilibrium and applying the notation of $x(t)=u(t)-u^*$, \eqref{VLm} is transformed to
\begin{align}\label{VLX}
 dx(t)  =X^{\ast}(t) \left\{\left[  -\tilde{A}x(t)-{\mathcal{A}}^d\tilde{{x}}^d(t)-{\mathcal{A}^D}\tilde{{x}}^D(t))\right]dt+g(x(t))dw(t)\right\},
\end{align}
where
 the notations
 \begin{align}
   X^{\ast}(t)&=\diag\{x_1(t)+u_1^{\ast}, \ldots,x_n(t)+u_n^{\ast}\}, \hspace{0.5cm}
 e = [ 1,  \ldots ,  1 ] ^T \in \mathbf{R}^{n \times 1}, \hspace{0.5cm} {\mathcal{I}}= e\otimes I_n
 \label{Jel-1} \\
  \tilde{x}^d(t)&=x^d(t)-{\mathcal{I}} x(t),\hspace{3.3cm}
   \tilde{x}^D(t)=x^D(t)-{\cal{B}} {\mathcal{I}} x(t) \label{Jel-2}
\end{align}
are applied, ${\mathcal{B}} \in\mathbf{R}^{n^2 \times n^2}$ is a diagonal matrix with
diagonal elements ${\mathcal{B}}_{kk}=\beta_{ij},$ for $ k=(i-1)n+j, \ i,j=1,\ldots ,n$
and the vectors $x^d(t), x^D(t) \in \mathbf{R}^{n^2 }$ are defined analogously to \eqref{Gjel1}.

In order to formulate the stability condition, let us define the block entry matrices
\begin{align*}\label{ek}
  {\mathcal{A}}=&[\tilde{A}, {\cal{A}}^d, {\cal{A}}^D, 0_{n \times n^2}]\in \mathbf{R}^{n \times (n+3n^2)}, \hspace{1.6cm}   
  &e_1=&[I_n, 0_{n\times n^2},0_{n\times n^2},0_{n\times n^2}]\in \mathbf{R}^{{n}\times (n+3n^2)},\\
  e_2=&[ 0_{n^2 \times n},I_{n^2},0_{n^2},0_{n^2}] \in \mathbf{R}^{{n^2}\times (n+3n^2)},
  &e_3=&[0_{n^2 \times n}, 0_{n^2},I_{n^2},0_{n^2}] \in \mathbf{R}^{{n^2}\times (n+3n^2)},\\
  e_4=&[0_{n^2 \times n}, 0_{n^2},0_{n^2},I_{n^2}]\in \mathbf{R}^{{n^2}\times (n+3n^2)},
  &e_5=&[{\mathcal{I}}, 0_{n^2},0_{n^2},0_{n^2}]\in \mathbf{R}^{{n^2}\times (n+3n^2)},
\end{align*}
and the diagonal matrices
$U^*=\diag\{u_1^*, \ldots, u_n^*\}\in\mathbf{R}^{n \times n}$, ${\mathcal{T}}, {\mathcal{T}}^d,  {\cal{A}}^{\alpha} \in\mathbf{R}^{n^2 \times n^2}$  with diagonal entries
 \begin{align*}
   {\mathcal{T}}_{kk}&=\overline{\tau}_{ij}, \hspace{0.5cm} {\mathcal{T}}_{kk}^d=1-\overline{\tau}_{ij}^d, \hspace{0.5cm} {\cal{A}}^{\alpha}_{kk}=\alpha_{ij}, \hspace{0.5cm}    k=(i-1)n+j, \ i,j=1,\ldots ,n.
 \end{align*}
\begin{thm} \label{thm:2}
Assume that the delay functions satisfy conditions \eqref{tau1} and \eqref{l-th1}. If there exist positive numbers $p_i, \; q_{ij}, \; r_{ij}, \; s_{ij}$ for $i,j=1, \ldots,n$ such that diagonal matrices
$P\in \mathbf{R}^{n \times n}$, ${\mathcal{Q}}, {\mathcal{R}}, {\mathcal{S}} \in\mathbf{R}^{n^2 \times n^2}$ with diagonal entries
\begin{equation}\label{th2-l1}
    P_{ii}=p_i,\hspace{0.5cm} {\mathcal{Q}}_{kk}=q_{ij}, \hspace{0.5cm}
    {\mathcal{R}}_{kk}=r_{ij}, \hspace{0.5cm}  {\mathcal{S}}_{kk}=s_{ij},  \hspace{0.5cm}  k=(i-1)n+j, \ i,j=1,\ldots ,n
\end{equation}
 satisfy the LMI
 \begin{equation}\label{th2-l2}
    \Sigma = \Sigma_1 + \Sigma_2 + \Sigma_3 + \Sigma_4 <0,
 \end{equation}
 where
 \begin{align}
    \Sigma_1=&
            -\frac{1}{2}e_1^T P {\cal{A}}-\frac{1}{2}{\cal{A}}^T P e_1+\frac{1}{2}e_1^T\sigma^T P U^* \sigma e_1 , \label{th2-l3}\\
    \Sigma_2=& e_1^T {\cal{Q}}_1 e_1- (e_2+e_5)^T {\cal{Q}}_d (e_2+e_5),  
    \hspace{1cm} {\cal{Q}}_1={\mathcal{I}}^T{\cal{Q}}{\mathcal{I}}, \hspace{0.5cm} {\cal{Q}}_d= {\mathcal{T}}^d {\cal{Q}},\label{th2-l4} \\
     \Sigma_3=& e_1^T {\cal{R}}_{1\ot} e_1-\begin{bmatrix}
                                         (e_3+{\cal{B}}e_5)^T, & e_4^T
                                       \end{bmatrix}
  \begin{bmatrix}
            4{\cal{R}}_{2\ot} & -6{\cal{R}}_{2\ot}  \\
                        -6{\cal{R}}_{2\ot} & 12{\cal{R}}_{2\ot}+4{\cal{A}}_{\alpha} {\cal{R}}
          \end{bmatrix}
          \begin{bmatrix}
            (e_3+{\cal{B}}e_5) \\
            e_4
          \end{bmatrix},   \label{th2-l5}\\
     \Sigma_4=&  e_1^T{\cal{S}}_1 e_4+e_4^T{\cal{S}}_{\clr{1}} e_1 -(e_2+{\cal{B}}e_5)^T {\cal{S}}_{2\ot}(e_2+{\cal{B}}e_5)
               -2e_4^T{\cal{A}}_{\alpha}{\cal{S}} e_4,  \label{th2-l6} \\
               & \hspace{0.45cm}{\cal{R}}_{1\ot}={\mathcal{I}}^T{\mathcal{T}}{\cal{R}}{\mathcal{I}}, \hspace{0.50cm}
          {\cal{R}}_{2\ot} = {\mathcal{T}}^{-1} {\cal{R}}, 
      \hspace{0.5cm} {\cal{S}}_1= {\mathcal{I}}^T{\cal{S}}, \hspace{0.5cm}  
      {\cal{S}}_{2\ot} ={\mathcal{T}}^{-1} {\cal{S}},  \label{thm:Sek}
 \end{align}
 then the equilibrium state $u^{\ast}$ of system \eqref{VLm} is almost surely globally asymptotically stable.
 \end{thm}
\textbf{Proof.}
Consider the Lyapunov-Krasovskii functional candidate 
\begin{align}
 &\hspace{1cm} V(t, x_t) = \sum_{\ell=1}^{4}V_\ell(t, x_t)= \sum_{i=1}^{n}V_1^i(x_i(t))+
   \sum_{\ell=2}^{4}\sum_{i,j=1}^{n} V_\ell^{ij}(t, x_t), \label{LKF} \\
  V_1^i(x_i)&=p_i\left(x_i-u_i^*\ln \frac{x_i+u_i^*}{u_i^*}\right), \hspace{0.6cm} 
  V_3^{ij}(t, x_t)= r_{ij} \int_{t-\overline{\tau}_{ij}}^{t}(\eta-t+\overline{\tau}_{ij})e^{2\alpha_{ij}(\eta-t)}x_j(\eta)^2 d\eta,
  \nonumber \\
  V_2^{ij}(t, x_t)&=q_{ij}\int_{t-{\tau}_{ij}(t)}^{t} x_j(\eta)^2d\eta,
  \hspace{1.1cm}
  V_4^{ij}(t, x_t)=\frac{ s_{ij}}{\overline{\tau}_{ij}}
  \int_{t-\overline{\tau}_{ij}}^{t}\left(\int_{\eta}^{t}e^{\alpha_{ij}(\zeta-t)}x_j(\zeta)d\zeta\right)^2 d\eta.  \nonumber
\end{align}
Fix an arbitrary $\varphi_0 \in C([-\overline{\tau},0],\mathbf{R}_{+}^{n})$ for \eqref{VLm}, and consider the solution  of \eqref{VLX} corresponding to the initial function $\phi_0(t)=\varphi_0(t)-u^{\ast},$
and  compute the derivative of $V$ applying the functional It\^{o}'s formula given in Equ. (3) of \cite{Zong2018}.
(Details and notations see in \cite{Zong2018}.)

Since {\clr{$V_1$ is independent of $t$, }} $\frac{\partial V_1}{\partial x_i} (x)=\frac{p_i x_i}{x_i+u_i^{\ast}}$, and $\frac{\partial^2 V_1}{\partial  x_i^2} (x)=\frac{p_i u_i^*}{(x_i(t)+u_i^*)^2}$, applying the above mentioned It\^{o}'s formula in symbolic differential form
 to $V_1(x)$ yields
\begin{align}\label{dV11}
 {\clr{ dV_1(x(t)) }}&=
  \left[-x(t)^T P \left(\tilde{A} x(t)+ {\cal{A}}^d \tilde{x}^d(t)+ {\cal{A}}^D \tilde{x}^D(t)\right)
  +   x(t)^T \sigma ^T P U^{\ast} \sigma x(t)  \right] dt + G(x(t)) dw(t),
\end{align}
where
$G(x)=[ p_1 x_1\Sigma_{j=1}^n \sigma_{1j}x_j,\ldots,
 p_n x_n\Sigma_{j=1}^n \sigma_{nj}x_j ] \in \mathbf{R}^{1\times n}. $

{\clr{Next one can verify that $\nabla_x V_\ell(t, x_t)=0,$ and $\mathcal{D} V_\ell(t, x_t)$ can be computed by taking the time derivative, if $\ell=2,3,4.$  Therefore $\mathcal{L}{V}_2(t,x_t)$ can be estimated as follows:}}
\begin{align}
{\clr{\mathcal{L}{V}_2(t,x_t)}}\leq & \sum_{j=1}^n \left(\sum_{i=1}^n q_{ij} \right)x_j(t)^2-
\sum_{i=1}^n \sum_{j=1}^n (1-{\ot}_{ij}^d)q_{ij}x_j(t-\tau_{ij}(t))^2 
 = x(t)^T {\mathcal{Q}}_1 x(t)- x^d(t)^T {\mathcal{Q}}_{d} x^d(t), \label{dv2ij}
\end{align}
where ${\mathcal{Q}}_1$ and ${\mathcal{Q}}_{d}$ are given in \eqref{th2-l4}, {\clr{and \eqref{tau1} has been taken into account.}}

Let us compute now ${\clr{\mathcal{L}{V}_3^{ij}(t,x_t)}}$:
\begin{equation}\label{dv3ij}
{\clr{\mathcal{L}{V}_3^{ij}(t,x_t)}}=r_{ij}\ot_{ij}x_j(t)^2-r_{ij}\int_{t-\ot_{ij}}^{t}e^{2\alpha_{ij}(\eta-t)}x_j(\eta)^2d\eta-2\alpha_{ij}r_{ij}\int_{t-\ot_{ij}}^{t}(\eta-t+\ot_{ij})e^{2\alpha_{ij}(\eta-t)}x_j(t)^2d\eta
\end{equation}
The first integral term can be estimated by the Wirtinger inequality (the applied form see in \cite{Lee2017}), while the second integral term can be estimated by the double{\clr{-}}integral Jensen inequality  (see e.g. \cite{Sun2009}). Thus, we obtain with $k=(i-1)n+j,$ $(i,j=1,\ldots,n)$ that
\begin{align}
  {\clr{\mathcal{L}{V}_3^{ij}(t,x_t)}} & \leq r_{ij}\ot_{ij}x_j(t)^2-\frac{r_{ij}}{\ot_{ij}}\left(({x_k^D}(t))^2+3({x_k^D}(t)-2{z_k^D}(t))^2\right)-4\alpha_{ij}r_{ij}({z_k^D}(t))^2 \nonumber \\&=
r_{ij}\ot_{ij}x_j(t)^2-\frac{r_{ij}}{\ot_{ij}}\left( 4({x_k^D}(t))^2-12{x_k^D}(t){z_k^D}(t)+12({z_k^D}(t))^2 \right) -4\alpha_{ij}r_{ij}({z_k^D}(t))^2. \label{dv3ije}
 \end{align}
where the vector
$z^D(t)\in \mathbf{R}^{n^2}$
is defined with the elements
\begin{equation}\label{zDij}
z_k^D(t)=\frac{1}{\overline{\tau}_{ij}}\int_{t-\overline{\tau}_{ij}}^{t}\int_{\eta}^{t}e^{\alpha_{ij}(\zeta-t)}x_j(\zeta)d\zeta d\eta . \hspace{0.5cm} 
\end{equation}
Summing up and applying \eqref{thm:Sek} we obtain that
\begin{align}\label{dV3M}
{\clr{\mathcal{L}{V}_3(t,x_t)}}&\leq x(t)^T{\cal{R}}_{1\ot}x(t)-\begin{bmatrix}
                                         x^D(t)^T, & z^D(t)^T
                                       \end{bmatrix}
  \begin{bmatrix}
            4{\cal{R}}_{2\ot} & -6{\cal{R}}_{2\ot}  \\
                        -6{\cal{R}}_{2\ot} & 12{\cal{R}}_{2\ot}+4{\cal{A}}_{\alpha}{\cal{R}}
          \end{bmatrix}
          \begin{bmatrix}
            x^D(t) \\
            z^D(t)
          \end{bmatrix} .
\end{align}

Computing  ${\clr{\mathcal{L}{V}_4^{ij}(t,x_t)}}$, and estimating the last term by Jensen's inequality yields
\begin{align}\label{dv4ij}
&{\clr{\mathcal{L}{V}_4^{ij}(t,x_t)}} 
=-\frac{s_{ij}}{\ot_{ij}}\left( \int_{t-\ot_{ij}}^{t}e^{\alpha_{ij}(\zeta-t)}x_j(\zeta)d\zeta\right)^2+
x_j(t)\frac{2s_{ij}}{\ot_{ij}}\int_{t-\ot_{ij}}^{t}\int_{\eta}^{t}e^{\alpha_{ij}(\eta-t)}x_j(\zeta)d\zeta d\eta \nonumber\\
& -
\frac{2\alpha_{ij}s_{ij}}{\ot_{ij}}\int_{t-\ot_{ij}}^{t}\left(\int_{\eta}^{t}e^{\alpha_{ij}(\zeta-t)}x_j(\zeta)d\zeta\right)^2 d\eta 
\leq 
-\frac{s_{ij}}{\ot_{uj}}(x^D_k(t))^2+2s_{ij}z^D_k(t)x_j(t)-2s_{ij}\alpha_{ij}(z^D_k(t))^2.
\end{align}
Summing up and applying \eqref{thm:Sek} we obtain that
\begin{align}\label{dV4M}
{\clr{\mathcal{L}{V}_4(t,x_t)}}&\leq -x^D(t)^T{\cal{S}}_{2\ot}x^D(t)+x(t)^T{\cal{S}}_1z^D(t)+z^D(t)^T{{\cal{S}}_1}x(t)-2z^D(t)^T{\cal{A}}_{\alpha}{\cal{S}}z^D(t). 
\end{align}
Let an extended variable $\xi(t)=(x(t)^T, \tilde{x}^d(t)^T, \tilde{x}^D(t)^T, {z}^D(t))^T\in \mathbf{R}^{n+3n^2}$ be introduced. Then, applying \eqref{dV11}-\eqref{dV4M} one can check with a straightforward computation that
\begin{equation}\label{thm2_veg}
  dV(t,x_t) \leq \xi(t)^T\left(\Sigma_1+\Sigma_2+\Sigma_3+\Sigma_4 \right)\xi(t)dt+G(x(t)) dw(t),
\end{equation}
where $\Sigma_{\ell}, \; \ell=1,\ldots,4$ are give by \eqref{th2-l3}-\eqref{thm:Sek}. Therefore, the statement of the theorem follows
from \eqref{th2-l2} based on the stability theory of stochastic differential equations (\cite{maokonyv}, \cite{mao2002}).
\begin{rem} Observe that it follows from \eqref{th2-l2} and \eqref{thm2_veg} that
   $\mathbb{E} V(T,x_T) \leq \mathbb{E} V(0,x_0)$ 
   for any $T>0.$
Therefore the conditions of Theorem \ref{thm:2} yield as an alternative for the existence of global positive solutions of \eqref{VLm}
for any initial function $\varphi_0 \in C([-\overline{\tau},0],\mathbf{R}_{+}^{n})$
with probability 1. (Cf. with Theorem 2.1 and Theorem 2.3 of \cite{mao2005}.) This alternative is useful in such cases, when the conditions $\sigma_{ii}>0, \; \forall i,$ are not satisfied (see Example 2 below).
\end{rem}
\section{Numerical examples}
\textbf{Example 1}. Consider a 3 species Lotka-Volterra model with the data
{\footnotesize
\begin{displaymath}
\begin{split}
A=\begin{bmatrix}2 & 1 & 0 \\ 0.5 & 2.5 & 0.5 \\ 0 & 1 & 2.5
 \end{bmatrix}, \hspace{1.0cm}
A^d=\begin{bmatrix}0.5 & 0.2 & 0.1 \\ 0.4 & 0.6 & 0 \\ 0.1 & 0 & 0.8
 \end{bmatrix}, \hspace{0.5cm}
A^D={\clr{\lambda_1}}\begin{bmatrix}0.4 & 0.5 & 0 \\ 0.2 & 1 & 0.1 \\ 0.1 & 0.1 & 0.5
 \end{bmatrix}, \\
 {\mathcal{T}}=\ot\begin{bmatrix}0.9 & 0.5 & 0.05 \\ 0.4 & 1 & 0.05 \\ 0.05 & 0.1 & 0.5
 \end{bmatrix}, \hspace{0.5cm}
  {\mathcal{T}}^d=\ot_d\begin{bmatrix}1 & 0.8 & 0.5 \\ 0.6 & 0.7 & 0.4 \\ 0.4 & 0.3 & 0.5
 \end{bmatrix}, \hspace{0.5cm}
 \sigma={\clr{\lambda_2}}\begin{bmatrix}0.2 & 0.05 & 0 \\ 0.15 & 0.1 & 0 \\ 0 & 0 & 0.2
 \end{bmatrix}
\end{split}
\end{displaymath}
}
and $\alpha_{ij}=2.$ The results obtained by {\clr{Theorem}} \ref{thm:2} for different values of the parameters ${\clr{\lambda_1}}, \; {\clr{\lambda_2}}, \; \ot $ and $\ot _d$ are given in Table 1.
\begin{table}[!ht] \label{Tab:num2}
\caption{{\footnotesize Maximum allowable delay bound $\overline{\tau}$ obtained by Theorem \ref{thm:2} }}
{\footnotesize
\begin{center}
\begin{tabular}{cc|ccccc}
\toprule
${\clr{\lambda_2}}$ & ${\clr{\lambda_1}}$ & $\ot_d=0$ & $\ot_d=0.2$ & $\ot_d=0.4$ & $\ot_d=0.6$ &   $\ot_d=0.6515$
\\
\hline
$1$ & $0$     & $100$    & $100$     & $100$     & $100$                        & $100$
            \\
    & $0.5$   & $3.5225$ & $2.4415$ & $1.5075$ & $0.4775$          &   $0.1725$
      \\
    & $1$     & $1.2745$ & $1.0065$ & $0.6675$ & $0.2305$        &     $0.0475$
      \\
\hline
$2$ & $0$     & $100$    & $100$    & $100$    & $100$           & $100$
            \\
    &$0.5$   &  $3.1075$ & $2.1575$ & $1.2775$ & $0.2950$        &   $0.0017$
      \\
    & $1$     & $1.1775$ & $0.9125$ & $0.5795$ & $0.1425$         &     infeasible
      \\
  \bottomrule
 \end{tabular}  \\
 \end{center}
  \vskip-2mm
     }
   \end{table}
We note that, in the case of $A^D=0,$ the LMI of Theorem \ref{thm:2} is independent of the delay upper bound, but it depends {\clr{only}} on the delay derivative upper bound. Thus, if it has a feasible solution for some $\ot_d$, $\ot, $
then it has a feasible solution for the same $\ot_d$ and arbitrary $\ot.$ {\clr{For the value $\ot_d=0.6515$ in the last column of Table 1, LMI \eqref{th2-l2} has a feasible solution, if $A^D=0$ and $\ot=100$, but it turns to be infeasible for any $\ot$, if $\ot_d$ is slightly increased.
 }}

\textbf{Example 2}. Consider the two species model of \cite{XiongAML2019}, where no delays are taken into account, i.e. $A^d=A^D=0,$ and let
{\footnotesize
\begin{equation*}
  A_1=\begin{bmatrix}3 & 1  \\ 2 & 2
 \end{bmatrix}, \hspace{0.3cm}
 \rho^{\clr{1}}=\begin{bmatrix}2  \\ 2
 \end{bmatrix}, \hspace{1.0cm}
 A_2=\begin{bmatrix}3 & 2  \\ 1 & 2
 \end{bmatrix}, \hspace{0.3cm}
 \rho^{\clr{2}}=\begin{bmatrix}5  \\ 1
 \end{bmatrix}, \hspace{1.0cm}
 \sigma=\begin{bmatrix}0 & \sigma_1  \\ \sigma_2 & 0
 \end{bmatrix}. \hspace{1.0cm}
\end{equation*}
}
If we take $A=A_1,$ $\rho=\rho^{\clr{1}},$
 $\sigma_1=1.5,$ $\sigma_2=2$, as in \cite{XiongAML2019}, and take formally the delay parameters as  ${\mathcal{T}}=10^{-5}I_2, $
  ${\mathcal{T}}^d=0, $ $\alpha_{\clr{ij}}=0,  $ then the LMI of Theorem \ref{thm:2} is feasible. If we change the data to $A=A_2, \rho=\rho^{\clr{2}},$
 $\sigma_1=\sqrt{2},$ $\sigma_2=\sqrt{2}$, then the LMI of Theorem \ref{thm:2} is feasible, while assumption $(H_1)$ of \cite{XiongAML2019} fails. Simulation supports the stability of the equilibrium $u^{\ast}=[1, \ 1]^T$.
 This suggest that Theorem \ref{thm:2} may lead to less conservative result, than some previously published stability conditions.

%

\section{Acknowledgements}
The research reported in this paper has been supported by the National Research, Development and Innovation Fund (TUDFO/51757/2019-ITM, Thematic Excellence Program).


\end{document}